\documentclass[12pt]{article}
\usepackage[tbtags]{amsmath}
\usepackage[]{amssymb}
\usepackage[]{amsbsy}
\usepackage[]{amsthm}
\usepackage[dvips]{graphicx}

\newtheorem{theorem}{Theorem}[section]

\newtheorem{lemma} {Lemma} [section]
\newtheorem{defi}{Definition} [section]
\newtheorem{example}{Example} [section]
\newtheorem{remark}{Remark} [section]

\def\<{\langle}
\def\>{\rangle}

\def\be{\begin{equation}}
\def\bea{\begin{eqnarray}}
\def\d{\delta}

\def\ee{\end{equation}}
\def\eea{\end{eqnarray}}

\def\splt#1{\begin{equation} \begin{split} #1 \end{split}\end{equation}}

\newcommand {\pa}{\partial}

\renewcommand{\(}{\left(}
\renewcommand{\)}{\right)}
\newcommand{\trho}{\tilde{\rho}}
\newcommand{\tiota}{\tilde{\iota}}

\begin{document}
\date{}
\baselineskip=1.2\normalbaselineskip
\thispagestyle{empty}

\title{{\noindent \bf  Boundary Feedback Control of Complex Ginzburg-Landau
Equation with A Simultaneously Space and Time Dependent
Coefficient}}

\author{
{\sc   Junji Jia  } \\
 \normalsize{ Department of Applied Mathematics, The University of Western Ontario }\\
  \normalsize{London, ON, N6A 5B7, Canada }\\
 \normalsize{\tt Email: jjia5@uwo.ca}
\vspace*{-2ex}}
\maketitle

\vspace*{2ex}
\begin{abstract}
Linearized complex Ginzburg-Landau equation models various physical
phenomena and the stability controls of them are important. In this
paper, we study the control of the LCGLE with a simultaneously space
and time dependent coefficient by transforming it into a complex
heat equation. It is shown that under certain conditions on the
coefficient functions $a_2(\tilde{x},\tilde{t})$, the exponential
stability of the system at any rate can be achieved by boundary
control based on the state feedback. The kernels are {\it
explicitly} calculated as series of approximation and shown to be
twice differentiable by using the {\it method of dominant}. Both the
exponential stabilities of the systems with Dirichlet and Neumann
boundary conditions are strictly proven.
\end{abstract}



\thispagestyle{empty}

\newpage
\section{Introduction}
\setcounter {equation} {0} \setcounter {theorem} {0} \setcounter
{lemma} {0} \hspace{15pt} In this paper, we use the boundary
feedback control law, which is based on the back stepping
methodology, to stabilize the linearized complex Ginzburg-Landau
equation (LCGLE) \be \frac{\partial u(\tilde{x},\tilde{t})}{\partial
\tilde{t}}=a_1\frac{\partial^2 u(\tilde{x},\tilde{t})}{\partial
\tilde{x}^2}+a_3(\tilde{x})\frac{\partial
u(\tilde{x},\tilde{t})}{\partial
\tilde{x}}+a_2(\tilde{x},\tilde{t})u(\tilde{x},\tilde{t})\ee where
$\tilde{x}\in(0,L),~\tilde{t}\in (0,T),~u:(0,L)\times(0,T)\to
\mathbb{C},~a_1\in\mathbb{C}$ and $\Re(a_1)>0,~a_3\in C^1([0,L];
\mathbb{C})$ and $a_2\in C^1([0,L]\times[0,T];\mathbb{C})$ is a
simultaneously space and time dependent coefficient function. The
Dirichlet boundary conditions of the system are \bea
u(0,t)&=&p(t),\\
u(L,t)&=&0\eea where $p(t)$ is the control input. This LCGLE models
various physical phenomena, such as the amplitude equation in
pattern formation \cite{cross} and the reaction diffusion of two
chemicals in one dimension \cite{cross, dl}, in all of which the
controls of the systems are important.

Without losing any generality, we set $a_3(\tilde{x})=0$ since this
term could always be eliminated by a gauge transform \be
u\to\tilde{u}=ue^{-1/2\int_0^{\tilde{x}}a_3(y)/a_1dy}\ee to obtain
the equation we will concentrate on hereafter \be \frac{\pa
\tilde{u}(\tilde{x},\tilde{t})}{\pa \tilde{t}}=a_1\frac{\pa^2
\tilde{u}(\tilde{x},\tilde{t})}{\pa
\tilde{x}^2}+a_2(\tilde{x},\tilde{t})\tilde{u}(\tilde{x},\tilde{t})
\label{eq1} .\ee Equation (\ref{eq1}) is equivalent to a complex
linear heat equation. The feedback control of this equation where
the $a_1$ and $a_2$ are real constants was first addressed by
Boskovic, Krstic and Liu in \cite{bkl} and the instability of the
systems was also shown there. Later the plant coefficient was
generalized to a space dependent case $a_2=a_2(\tilde{x})$ by Liu
\cite{liunew}. More recently, Aamo, Smyshlyaev and Krstic
\cite{aamo1, kris2, kris3} considered the generalization to a
complex equation with a complex space or time dependent coefficient
function, i.e., $a_2=a_2(\tilde{x})$ or $a_2=a_2(\tilde{t})$. In
this paper, we extend their work by generalizing the plant
coefficient to a simultaneously space and time dependent case, i.e.,
$a_2=a_2(\tilde{x},\tilde{t})$. We also considered the Neumann
boundary condition problems. Besides, from the proof of lemma
(\ref{lm:kers}) we can see that the explicit result for the control
kernel is given and therefore it could be immediately used in
numerical simulation.

The paper is organized as follows. In section \ref{sec2}, we prove
that for some $a_2(\tilde{x},\tilde{t})$, the system (\ref{eq1})
with the Dirichlet boundary conditions could be stabilized by the
boundary feedback control. The unique existence of the kernels are
shown in subsection \ref{ssec21}; Using these kernels, in subsection
\ref{ssec22} the original system is transformed into a new
well-posed system. In subsection \ref{ssec23} we show that the
original system is exponentially stabilized. In section \ref{sec3},
we prove that the Neumann boundary problem can also be stabilized
using the similar methodology.

\section{Dirichlet Boundary Value Problem}\label{sec2}
\setcounter {equation} {0} \setcounter {theorem} {0} \setcounter
{lemma} {0} First we break equation (\ref{eq1}) into two coupled
PDEs with real domains, real coefficient functions and real
codomains by defining \bea
\rho(x,t)&=&\Re(u(x,t))=\frac{1}{2}\(u(x,t)+\overline{u}(x,t)\),\\
\iota(x,t)&=&\Im(u(x,t))=\frac{1}{2i}\(u(x,t)-\overline{u}(x,t)\)\eea
where \be x=\frac{L-\tilde{x}}{L}, t=\frac{\tilde{t}}{T} \mbox{ and
} u(x,t)=\tilde{u}(\tilde{x},\tilde{t})\label{eq:trans}\ee and $^-$
denotes the complex conjugate. Equation (\ref{eq1}) then becomes
\bea
\rho_t&=&a_R\rho_{xx}+b_R(x,t)\rho-a_I\iota_{xx}-b_I(x,t)\iota,\label{eq:sysrho}\\
\iota_t&=&a_R\iota_{xx}+b_R(x,t)\iota+a_I\rho_{xx}+b_I(x,t)\rho\label{eq:sysiota}\eea
for $(x,t)\in(0,1)\times (0,1)\equiv \Sigma$, with boundary
conditions
\bea \rho(0,t)&=&0,~\iota(0,t)=0\label{eq:lbc}\\
\rho(1,t)&=&p_R(t),~\iota(1,t)=p_I(t)\label{eq:rbc}\eea
where \bea
a_R&=&\frac{1}{L^2}\Re{(a_1)},~a_I=\frac{1}{L^2}\Im{(a_1)},\label{eq:arcon}\\
b_R&=&\Re{\(a_2(\tilde{x},\tilde{t})\)},~b_I=\Im{\(a_2(\tilde{x},\tilde{t})\)}.\label{eq:bricon}\eea
Note that the transformation (\ref{eq:trans}) is to normalize the
$x$ domain and the $t$ domain. To stabilize the new equation system
(\ref{eq:sysrho})-(\ref{eq:bricon}), we set the boundary control
input $p(t)$ into the form \bea
p_R(t)&=&\int_0^1[k(1,y,t)\rho(y,t)+k_c(1,y,t)\iota(y,t)]dy\label{eq:inpr}\\
p_I(t)&=&\int_0^1[k(1,y,t)\iota(y,t)-k_c(1,y,t)\rho(y,t)]dy,\label{eq:inpi}\eea
where $k$ and $k_c$ are kernels that we should find.

\subsection{The Unique Existence of the Kernels}\label{ssec21}
First, we state a lemma about the kernels $k$ and $k_c$ for future
use. It will be clear in lemma (\ref{lm:coninv}) why we consider
this PDE system.
\begin{lemma}\label{lm:kers}
If $a_Rb_R(x,t)+a_Ib_I(x,t)$ is analytic in $t$, the partial differential
equation system of $k$ and $k_c$:
\bea
k_{xx}&=&k_{yy}+\beta(x,y,t)k+\beta_c(x,y,t)k_c+p_1k_t+p_2k_{c,t}\label{eq:kerk}\\
k_{c,xx}&=&k_{c,yy}-\beta_c(x,y,t)k+\beta(x,y,t)k_c+q_1k_{c,t}+q_2k_t\label{eq:kerkc}\eea
for $(x,y)\in \Omega$ with boundary conditions
\bea
k(x,x,t)&=&-\frac{1}{2}\int_0^x\beta(\gamma,\gamma,t)d\gamma\label{eq:kintcon}\\
k_c(x,x,t)&=&\frac{1}{2}\int_0^x\beta_c(\gamma,\gamma,t)d\gamma\label{eq:kcintcon}\\
k(x,0,t)&=&0\label{eq:kcond}\\
k_c(x,0,t)&=&0\label{eq:kccond}\eea
where
\bea
p_1&=&q_1=a_R/(a_R^2+a_I^2),\label{eq:pq1con}\\
p_2&=&-q_2=-a_I/(a_R^2+a_I^2),\label{eq:pq2con}\\
\beta(x,y,t)&=&[a_R(b_R(y,t)-f_R(x,t))+a_I(b_I(y,t)-f_I(x,t))]/(a_R^2+a_I^2),\label{eq:betacon}\\
\beta_c(x,y,t)&=&[a_R(b_I(y,t)-f_I(x,t))-a_I(b_R(y,t)-f_R(x,t))]/(a_R^2+a_I^2),\label{eq:betaccon}
\eea
in which $f_R(x,t),~f_I(y,t)$
defined in theorem (\ref{thm:dic}) are analytic in
$t$ and differentiable in $y$, has a unique solution satisfying
\bea
|k(x,y,t)|&\leq& Me^{M(x^2-y^2)}\label{eq:bdk}\\
|k_c(x,y,t)|&\leq& Me^{M(x^2-y^2)}\label{eq:bdkc}~,\eea where $M$ is a positive
constant, for any time interval $(0,t_0)(t_0<1)$ we are concerning.
\end{lemma}

\begin{proof}
Using the substitutions \bea
&x=\xi+\eta,~y=\xi-\eta,&\label{eq:xietas}\\
&G(\xi,\eta,t)=k(x,y,t),~G_c(\xi,\eta,t)=k_c(x,y,t),&\label{eq:gs}\\
&\delta(\xi,\eta,t)=\beta(x,y,t),~\delta_c(\xi,\eta,t)=\beta_c(x,y,t),&\label{eq:deltas}\eea
we can have \be \begin{array}{ll}
k_{xx}=1/4(G_{\xi\xi}+2G_{\xi\eta}+G_{\eta\eta})&\quad k_{c,xx}=1/4(G_{c,\xi\xi}+2G_{c,\xi\eta}+G_{c,\eta\eta})\\
k_{yy}=1/4(G_{\xi\xi}-2G_{\xi\eta}+G_{\eta\eta})&\quad k_{c,yy}=1/4(G_{c,\xi\xi}-2G_{c,\xi\eta}+G_{c,\eta\eta})\\
k_t=G_t&\quad k_{c,t}=G_{c,t}.\end{array}\ee Thus the system
(\ref{eq:kerk})-(\ref{eq:kccond}) is transformed into \bea
G_{\xi\eta}(\xi,\eta,t)&=&\(\delta(\xi,\eta,t)+p_1\frac{\pa}{\pa
t}\)G(\xi,\eta,t)+\(\delta_c(\xi,\eta,t)+p_2\frac{\pa}{\pa
t}\)G_c(\xi,\eta,t)\label{eq:intg}\\
G_{c,\xi\eta}(\xi,\eta,t)&=&\(\delta(\xi,\eta,t)+q_1\frac{\pa}{\pa
t}\)G_c(\xi,\eta,t)+\(-\delta_c(\xi,\eta,t)+q_2\frac{\pa}{\pa
t}\)G(\xi,\eta,t)\label{eq:intgc}\eea
with boundary conditions
\bea
G(\xi,0,t)&=&-\frac{1}{2}\int_0^{\xi}\delta(\tau,0,t)\d\tau,\label{eq:zerogc}\\
G_c(\xi,0,t)&=&\frac{1}{2}\int_0^{\xi}\delta_c(\tau,0,t)\d\tau,\label{eq:zerogcc}\\
G(\xi,\xi,t)&=&0,\label{eq:impbc}\\
G_c(\xi,\xi,t)&=&0.\eea By direct integrating of (\ref{eq:intg}) and
(\ref{eq:intgc}) with respect to $\xi$ and $\eta$ and using
(\ref{eq:zerogc})-(\ref{eq:zerogcc}), we obtain the equivalent
integral equations \splt{\label{eq:geq}
G(\xi,\eta,t)=-\frac{1}{2}\int_{\eta}^{\xi}\delta(\tau,0,t)d\tau&+\int_{\eta}^{\xi}\!\!\int_0^{\eta}\(\delta(\tau,s,t)+p_1\frac{\pa
}{\pa t}\)G(\tau,s,t)dsd\tau\\
&+\int_{\eta}^{\xi}\!\!\int_0^{\eta}\(\delta_c(\tau,s,t)+p_2\frac{\pa
}{\pa t}\)G_c(\tau,s,t)dsd\tau}
\splt{\label{eq:gceq}
G_c(\xi,\eta,t)=\frac{1}{2}\int_{\eta}^{\xi}\delta(\tau,0,t)d\tau&+\int_{\eta}^{\xi}\!\!\int_0^{\eta}\(\delta(\tau,s,t)+q_1\frac{\pa
}{\pa t}\)G_c(\tau,s,t)dsd\tau\\
&+\int_{\eta}^{\xi}\!\!\int_0^{\eta}\(-\delta_c(\tau,s,t)+q_2\frac{\pa}{\pa t}\)G(\tau,s,t)dsd\tau}
\\
To proceed, we set \splt{\label{eq:split}
G_0(\xi,\eta,t)&=-\frac{1}{2}\int_{\eta}^{\xi}\delta(\tau,0,t)\d\tau,~G_{c,0}(\xi,\eta,t)=\frac{1}{2}\int_{\eta}^{\xi}\delta(\tau,0,t)\d\tau}
\splt{\label{eq:gn1eq}
G_{n+1}(\xi,\eta,t)&=\int_{\eta}^{\xi}\!\!\int_0^{\eta}\(\delta(\tau,s,t)+p_1\frac{\pa}{\pa
t}\)G_n(\tau,s,t)dsd\tau \\
&+\int_{\eta}^{\xi}\!\!\int_0^{\eta}\(\delta_c(\tau,s,t)+p_2\frac{\pa}{\pa
t}\)G_{c,n}(\tau,s,t)dsd\tau}
\splt{\label{eq:gcn1eq}
G_{c,n+1}(\xi,\eta,t)&=\int_{\eta}^{\xi}\!\!\int_0^{\eta}\(\delta(\tau,s,t)+q_1\frac{\pa}{\pa
t}\)G_{c,n}(\tau,s,t)dsd\tau\\
&+\int_{\eta}^{\xi}\!\!\int_0^{\eta}\(-\delta_c(\tau,s,t)+q_2\frac{\pa}{\pa
t}\)G_n(\tau,s,t)dsd\tau.}
We want to show that both the series
\be G=\sum_{n=1}^{\infty}G_n(\xi,\eta,t)\mbox{ and
}G_c=\sum_{n=1}^{\infty}G_{c,n}(\xi,\eta,t)\ee
are uniquely convergent and thus by the method of successive approximation they are the solutions of system
(\ref{eq:geq})-(\ref{eq:gceq}).

This could be done by the virtue of Colton's {\it method of
dominant} \cite{col1, col2}. This method works as follows. If we are
given two series \be
S_1=\sum_{n=1}^{\infty}a_{1n}t^n,~S_2=\sum_{n=1}^{\infty}a_{2n}t^n,~t\in
(0,1),\ee where $a_{2n}\geq0$, then we say $S_2$ {\it dominates}
$S_1$ if $|a_{1n}|\leq a_{2n},~n=1,2,3,\cdots$, and write $S_1\ll
S_2$. It can be easily checked that \bea
\mbox{if } S_1\ll S_2,&&\mbox{then }|S_1|\leq S_2,\label{eq:abspro}\\
&&\mbox{ and }\frac{\pa S_1}{\pa t}\ll\frac{\pa S_2}{\pa t},~S_1\ll S_2(1-t)^{-1};\\
\mbox{if } S_1\ll S_2,S_2\ll S_3,&&\mbox{then }S_1\ll S_3;\\
\mbox{if } S_1\ll S_2,S_3\ll S_4,&&\mbox{then }S_1+S_2\ll
S_3+S_4.\eea Moreover, the property of {\it dominancy} can also be
kept if the integrals of the two series are not with respect to $t$
but other variables, that is \be \mbox{if } S_1\ll S_2,\mbox{ then
}\int_a^bS_1dx\ll\int_a^bS_2dx.\label{eq:intpro}\ee On the other
hand, if a function $f$ is analytic in $t\in(0,1)$, then there exist
a positive constant $C$ such that $f\ll C(1-t)^{-1}$.

Using this method in our problem, it can be shown that equations (\ref{eq:kerk})
and (\ref{eq:kerkc}) have unique twice continuously differentiable solutions if
$\delta$ and $\delta_c$ are analytic in $t$. In fact, since
$a_Rb_R(x,t)+a_Ib_I(x,t)$ is analytic in $t$ and $f_I(x,t)$ and $f_R(x,t)$ can also be chosen to be analytic in $t$, from equation
(\ref{eq:betacon}), (\ref{eq:betaccon}) and
(\ref{eq:deltas}), we know that $\delta$ and $\delta_c$ are analytic in $t$.
Thus we can let $C,~C_c$ be two positive constants larger than 1 such that
\be \delta(\xi,\eta,t)\ll C(1-t)^{-1},~\delta_c(\xi,\eta,t)\ll C_c(1-t)^{-1},\ee
or further
\be \delta(\xi,\eta,t)\ll N(1-t)^{-1},~\delta_c(\xi,\eta,t)\ll N(1-t)^{-1},~N=\max\{C,C_c\}.\ee
From equation of (\ref{eq:split}) and noticing (\ref{eq:intpro}), we have \bea
G_0&\ll& N(1-t)^{-1}\\
G_{c,0}&\ll& N(1-t)^{-1}.\eea
By induction, suppose
\bea
G_n&\ll& \frac{4^n(\xi\eta)^nN^{n+1}}{n!}(1-t)^{-(n+1)}\\
G_{c,n}&\ll& \frac{4^n(\xi\eta)^nN^{n+1}}{n!}(1-t)^{-(n+1)},\eea
we will have from (\ref{eq:gn1eq})
\splt{
G_{n+1}=&\int_{\eta}^{\xi}\!\!\int_0^{\eta}\(\delta+p_1\frac{\pa}{\pa
t}\)G_ndsd\tau+\int_{\eta}^{\xi}\!\!\int_0^{\eta}\(\delta_c+p_2\frac{\pa}{\pa
t}\)G_{c,n}dsd\tau\\
\ll&\int_{\eta}^{\xi}\!\!\int_0^{\eta}\(\frac{N}{1-t}+p_1\frac{\pa }{\pa
t}\)\frac{4^n(s\tau)^nN^{n+1}}{n!}(1-t)^{-(n+1)}dsd\tau\\
&+\int_{\eta}^{\xi}\!\!\int_0^{\eta}\(\frac{N}{1-t}+p_2\frac{\pa }{\pa
t}\)\frac{4^n(s\tau)^nN^{n+1}}{n!}(1-t)^{-(n+1)}dsd\tau\\
\ll&\frac{4^n(\xi\eta)^{n+1}N^{n+1}}{((n+1)!)^2}(1-t)^{-(n+2)}(2N+np_1+np_2)\\
\ll&\frac{4^{n+1}(\xi\eta)^{n+1}N^{n+2}}{(n+1)!}(1-t)^{-(n+2)},} where all the properties
(\ref{eq:abspro})-(\ref{eq:intpro}) have been used.
Similarly, we also have
\be G_{c,n+1}\ll\frac{4^{n+1}(\xi\eta)^{n+1}N^{n+2}}{(n+1)!}(1-t)^{-(n+2)}\ee
and hence by (\ref{eq:abspro}) \be
|G_{n+1}|\leq\frac{4^{n+1}(\xi\eta)^{n+1}N^{n+2}}{(n+1)!}(1-t)^{-(n+2)}~,
|G_{c,n+1}|\leq \frac{4^{n+1}(\xi\eta)^{n+1}N^{n+2}}{(n+1)!}(1-t)^{-(n+2)}.\ee

It is clear that the two series convergent absolutely and uniformly and are
solution of (\ref{eq:geq})-(\ref{eq:gceq}). $G$ and $G_c$ are $C^2$ since
$\delta(\xi,\eta,t)$ and $\delta_c(\xi,\eta,t)$ are $C^1$. Since
\be 4^n(\xi\eta)^nN^n(1-t)^{-n}=\((x^2-y^2)\frac{N}{1-t}\)^n,\label{eq:bornb}\ee
(\ref{eq:bdk}) and (\ref{eq:bdkc}) follow directly by assigning
$N/(1-t_0)=M$.
\end{proof}
\begin{remark}
(1). The proof of Lemma (\ref{lm:kers}) provides a numerical
computation scheme of successive approximation to compute the
kernels. This makes the feedback laws very useful in real problems.
(2). We require $t\in (0,t_0),~t_0<1$ in (\ref{eq:bdk}) and
(\ref{eq:bdkc}) since from (\ref{eq:bornb}) one can see that $k$ and
$k_c$ are not bounded when $t\to1$. This requirement is tolerable
since in practice, we only require the system to be stable in a time
interval $\tilde{t}\in(0,T_0)$ where $T_0=t_0T$ from rescaling
equation (\ref{eq:trans}) and the $T$ can be choose as large as we
want.
\end{remark}

\subsection{Conversion and Inverse Convertibility of the Original
System}\label{ssec22}
We want to show that the original system (\ref{eq:sysrho})-(\ref{eq:inpi}) can be converted to a new system
by integral transform with $k$ and $k_c$ as kernels and this conversion is invertible.
\begin{lemma}\label{lm:coninv}
Let $k(x,y,t)$ and $k_c(x,y,t)$ be the solution of
(\ref{eq:kerk})-(\ref{eq:kerkc}) and define a pair of linear bounded operator
$K$ and $K_c$ by
\bea
\trho(x,t)&=&(K\rho)(x,t)=\rho(x,t)-\int_0^x[k(x,y,t)\rho(y,t)+k_c(x,y,t)\iota(y,t)]dy\label{eq:rhotrs}\\
\tiota(x,t)&=&(K_c\iota)(x,t)=\iota(x,t)-\int_0^x[-k_c(x,y,t)\rho(y,t)+k(x,y,t)\iota(y,t)]dy,\label{eq:iotatrs}\eea
Then,\\
1.$K$ and $K_c$ convert system (\ref{eq:sysrho})-(\ref{eq:inpi}) to system
\bea
\trho_t&=&a_R\trho_{xx}+f_R(x,t)\trho-a_I\tiota_{xx}-f_I(x,t)\tiota,\label{eq:trrho}\\
\tiota_t&=&a_I\trho_{xx}+f_I(x,t)\trho+a_R\tiota_{xx}+f_R(x,t)\tiota,\label{eq:triota}\eea
for $x\in (0,1)$, with boundary conditions
\be \trho(0,t)=0,~\tiota(0,t)=0,~\trho(1,t)=0,~\tiota(1,t)=0.\label{eq:bcnew}\ee
2. Both $K$ and $K_c$ have linear bounded inverses.
\end{lemma}

\begin{proof}
{\it Proof of part 1.}: To prove {\it 1}, we compute as follows.
Differentiating (\ref{eq:rhotrs}) and using
(\ref{eq:sysrho})-(\ref{eq:sysiota}), we have \splt{\nonumber
\trho_t(x,t)=&a_R\rho_{xx}+b_R(x,t)\rho-a_I\iota_{xx}-b_I(x,t)\iota\\
&-\int_0^x[k(x,y,t)(a_R\rho_{xx}+b_R(y,t)\rho-a_I\iota_{xx}-b_I(y,t)\iota)+k_t\rho\\
&\quad\quad
k_c(x,y,t)(a_I\rho_{xx}+b_I(y,t)\rho+a_R\iota_{xx}+b_R(y,t)\iota)+k_{c,t}\iota]dy.}
Using (\ref{eq:rhotrs})-(\ref{eq:iotatrs}), integrating the integral
by parts twice and using (\ref{eq:lbc}),
(\ref{eq:kcond})-(\ref{eq:kccond}) and (\ref{eq:bcnew}), we will get
\splt{\nonumber \trho_t(x,t)=&a_R\(\trho_{xx}(x,t)+\frac{\pa^2}{\pa
x^2}\int_0^x[k(x,y,t)\rho(y,t)+k_c(x,y,t)\iota(y,t)]dy\)\\
+&b_R(x,t)\(\trho(x,t)+\int_0^x[k(x,y,t)\rho(y,t)+k_c(x,y,t)\iota(y,t)]dy\)\\
-&a_I\(\tiota_{xx}(x,t)+\frac{\pa^2}{\pa
x^2}\int_0^x[-k_c(x,y,t)\rho(y,t)+k(x,y,t)\iota(y,t)]dy\)\\
-&b_I(x,t)\(\tiota(x,t)+\int_0^x[-k_c(x,y,t)\rho(y,t)+k(x,y,t)\iota(y,t)]dy\)}
\splt{
&-k(x,x,t)a_R\rho_x(x,t)+k(x,x,t)a_I\iota_x(x,t)-k_c(x,x,t)a_I\rho_x(x,t)-k_c(x,x,t)a_R\iota_x(x,t)\\
&+k_y(x,x,t)a_R\rho(x,t)-k_y(x,x,t)a_I\iota(x,t)+k_{c,y}(x,x,t)a_I\rho(x,t)+k_{c,y}(x,x,t)a_R\iota(x,t)\\
&-\int_0^x[k_{yy}(x,y,t)(a_R\rho(y,t)-a_I\iota(y,t))+k(x,y,t)(b_R(y,t)\rho(y,t)-b_I(y,t)\iota(y,t))\\
&\quad\quad+k_{c,yy}(x,y,t)(a_I\rho(y,t)+a_R\iota(y,t))+k_c(x,y,t)(b_I(y,t)\rho(y,t)+b_R(y,t)\iota(y,t))\\
&\quad\quad+k_t(x,y,t)\rho(y,t)+k_{c,t}(x,y,t)\iota(y,t)]dy.}
Using the identity
\splt{
\frac{\pa^2}{\pa
x^2}\int_0^x\lambda(x,y,t)\chi(y,t)dy=&\int_0^x\lambda_{xx}(x,y,t)\chi(y,t)dy+\lambda_x(x,x,t)\chi(x,t)\\
&+\chi(x,t)\frac{\pa \lambda(x,x,t)}{\pa x}+\lambda(x,x,t)\chi_x(x,t)} and the
equation (\ref{eq:trrho})-(\ref{eq:triota}) again, we obtain
\splt{\label{eq:trhofn}
\trho_t(x,t)=&a_R\trho_{xx}(x,t)-a_I\tiota_{xx}(x,t)+b_R(x,t)\trho(x,t)-b_I(x,t)\tiota(x,t)\\
&+2\(a_R\frac{\pa k(x,x,t)}{\pa x}+a_I\frac{\pa k_c(x,x,t)}{\pa
x}\)\trho(x,t)+\int_0^xR(x,y,t)\rho(y,t)dy\\
&+2\(a_R\frac{\pa k_c(x,x,t)}{\pa x}-a_I\frac{\pa k(x,x,t)}{\pa
x}\)\tiota(x,t)+\int_0^xI(x,y,t)\iota(y,t)dy}
where
\splt{\label{eq:recof}
R(x,y,t)=&a_R(k_{xx}(x,y,t)-k_{yy}(x,y,t))+a_I(k_{c,xx}(x,y,t)-k_{c,yy}(x,y,t))-k_t(x,y,t)\\
&+2\(a_R\frac{\pa k(x,x,t)}{\pa x}+a_I\frac{\pa k_c(x,x,t)}{\pa
x}+b_R(x,t)-b_R(y,t)\)k(x,y,t)\\
&+2\(-a_R\frac{\pa k_c(x,x,t)}{\pa x}+a_I\frac{\pa k(x,x,t)}{\pa
x}+b_I(x,t)-b_I(y,t)\)k_c(x,y,t),}
\splt{\label{eq:imcof}
I(x,y,t)=&a_R(k_{c,xx}(x,y,t)-k_{c,yy}(x,y,t))-a_I(k_{xx}(x,y,t)-k_{yy}(x,y,t))-k_{c,t}(x,y,t)\\
&+2\(a_R\frac{\pa k(x,x,t)}{\pa x}+a_I\frac{\pa k_c(x,x,t)}{\pa
x}+b_R(x,t)-b_R(y,t)\)k_c(x,y,t)\\
&+2\(a_R\frac{\pa k_c(x,x,t)}{\pa x}-a_I\frac{\pa k(x,x,t)}{\pa
x}-b_I(x,t)+b_I(y,t)\)k(x,y,t).}
Substituting equation (\ref{eq:kerk})-(\ref{eq:kerkc}) and (\ref{eq:kintcon})-(\ref{eq:kcintcon}) into
(\ref{eq:recof})-(\ref{eq:imcof}), we can get
\splt{
&R(x,y,t)=[(a_Rp_1+a_Iq_2)-1]k_t(x,y,t)+(a_Rp_2+a_Iq_1)k_{c,t}(x,y,t))\\
&+(a_R\beta(x,y,t)-a_I\beta_c(x,y,t)-(a_R\beta(x,x,t)-a_I\beta_c(x,x,t))+b_R(x,t)-b_R(y,t))k(x,y,t)\\
&+(a_R\beta_c(x,y,t)+a_I\beta_c(x,y,t)-(a_R\beta_c(x,x,t)+a_I\beta(x,x,t))+b_I(x,t)-b_I(y,t))k_c(x,y,t),}
\splt{
&I(x,y,t)=[(a_Rq_1-a_Ip_2)-1]k_{c,t}(x,y,t)+(a_Rq_2+a_Ip_1)k_t(x,y,t))\\
&+(-a_R\beta_c(x,y,t)-a_I\beta(x,y,t)+(a_R\beta_c(x,x,t)+a_I\beta(x,x,t))-b_I(x,t)+b_I(y,t))k(x,y,t)\\
&+(a_R\beta(x,y,t)-a_I\beta_c(x,y,t)-(a_R\beta(x,x,t)-a_I\beta_c(x,x,t))+b_R(x,t)-b_R(y,t))k_c(x,y,t).}
From (\ref{eq:betacon})-(\ref{eq:betaccon}), we have
\bea
a_R\beta(x,y,t)-a_I\beta_c(x,y,t)&=&b_R(y,t)-f_R(x,t)\label{eq:fromb1}\\
a_R\beta_c(x,y,t)+a_I\beta(x,y,t)&=&b_I(y,t)-f_I(x,t)\label{eq:fromb2}\eea
and with the help of (\ref{eq:pq1con})-(\ref{eq:pq2con}), it is easy to
check that
\be R(x,y,t)=I(x,y,t)\equiv 0.\ee
Thus equation (\ref{eq:trhofn}) becomes the following after using
(\ref{eq:kintcon}), (\ref{eq:kcintcon}) again,
\splt{
\trho(x,t)=&a_R\trho_{xx}(x,t)-a_I\tiota_{xx}(x,t)\\
&+(-a_R\beta(x,x,t)+a_I\beta_c(x,x,t)+b_R(x,t))\trho(x,t)\\
&+(a_R\beta_c(x,x,t)+a_I\beta(x,x,t)-b_I(x,t))\tiota(x,t).}
Substituting (\ref{eq:fromb1}) and (\ref{eq:fromb2}) into above, we can obtain
(\ref{eq:trrho}). Equation (\ref{eq:triota}) follows similarly by
differentiating (\ref{eq:iotatrs}) and applying the same procedure. The boundary conditions (\ref{eq:bcnew}) follow from setting
$x=0$ and $x=1$ in (\ref{eq:rhotrs}) and (\ref{eq:iotatrs}) and using
(\ref{eq:lbc})-(\ref{eq:rbc}) and (\ref{eq:inpr})-(\ref{eq:inpi}). Thus we
proved part {\it 1}.

{\it Proof of part 2.}: To prove {\it 2}, the invertibility of the
transform, consider the following transform \bea
\rho(x,t)&=&(K^{-1}\trho)(x,t)=\trho(x,t)-\int_0^x[l(x,y,t)\trho(y,t)+l_c(x,y,t)\tiota(y,t)]dy\label{eq:bktrrho}\\
\iota(x,t)&=&(K_c^{-1}\tiota)(x,t)=\tiota(x,t)-\int_0^x[-l_c(x,y,t)\trho(y,t)+l(x,y,t)\tiota(y,t)]dy,\label{eq:bktriota}.\eea
They are actually the inverse transform if there exist unique solutions
for $l(x,y,t)$ and $l_c(x,y,t)$ that satisfying
\bea
l_{xx}&=&k_{yy}-\beta(y,x,t)k-\beta_c(y,x,t)k_c-p_1k_t-p_2k_{c,t}\label{eq:kerl}\\
l_{c,xx}&=&k_{c,yy}+\beta_c(y,x,t)k-\beta(y,x,t)k_c-q_1k_{c,t}-q_2k_t\label{eq:kerlc}\eea
with boundary conditions
\bea
l(x,x,t)&=&\frac{1}{2}\int_0^x\beta(\gamma,\gamma,t)d\gamma\label{eq:lintcon}\\
l_c(x,x,t)&=&-\frac{1}{2}\int_0^x\beta_c(\gamma,\gamma,t)d\gamma\label{eq:lcintcon}\\
l(x,0,t)&=&l_c(x,0,t)=0,\label{eq:lcrcon} \eea and \bea
|l(x,y,t)|&\leq& Me^{M(x^2-y^2)},\label{eq:bdl}\\
|l_c(x,y,t)|&\leq& Me^{M(x^2-y^2)}\label{eq:bdlc}~,\eea where $M$ is
a positive constant, for any time interval $(0,t_0) (t_0<1)$ we are
concerning.

The proof that the transform (\ref{eq:bktrrho}) and
(\ref{eq:bktriota}) will lead to system
(\ref{eq:kerl})-(\ref{eq:lcrcon}) is similar to the proof of part
{\it 1} of this lemma. It is also clear that the transforms are
linear. The unique existence of the solution to system
(\ref{eq:kerl})-(\ref{eq:bdlc}) is similar to the proof of lemma
(\ref{lm:kers}) and thus the transforms are also bounded. Thus this
proved that there exist the bounded linear inverse transforms as
shown in (\ref{eq:bktrrho}) and (\ref{eq:bktriota}).
\end{proof}

\subsection{Stability of the Controlled System}\label{ssec23}
Our purpose of all doing this is to control the original system by
the boundary inputs. Its stability is established in the following
theorem.

\begin{theorem}\label{thm:dic}
There exist feedback kernels $k(1,\cdot,t),~k_c(1,\cdot,t)\in
C^{2,\infty}(0,1)\times(0,t_0), t_0<1$, such that for arbitrary
initial data $\rho^0(x),~\iota^0(x)\in H^1(0,1)$, system
(\ref{eq:sysrho})-(\ref{eq:rbc}) with
(\ref{eq:inpi})-(\ref{eq:inpr}), where $b_R(x,t)$ and $b_I(x,t)$ are
analytic in $t$, has a unique solution $\rho,~\iota\in
C^{2,\infty}(0,1)\times(0,t_0)$ that is exponentially stable in the
$L_2(0,1)$ and $H_1(0,1)$ norms.
\end{theorem}
\begin{proof}
First we notice that problem (\ref{eq:sysrho})-(\ref{eq:rbc}) with
(\ref{eq:inpi})-(\ref{eq:inpr}) is well posed since, by lemma
(\ref{lm:coninv}) they can be transformed into the problem
(\ref{eq:trrho})-(\ref{eq:bcnew}) via the isomorphism defined by
(\ref{eq:rhotrs}) and (\ref{eq:iotatrs}) and problem
(\ref{eq:trrho})-(\ref{eq:bcnew}) is well posed (see, e.g.
\cite{wellp}).

Since through part {\it 2} of lemma (\ref{lm:coninv}), we showed
that there exist linear bounded invertible transforms, it is
sufficient to show that $(\trho,\tiota)$ are exponentially stable.
The proof of this could be done by the same definitions and
evaluations of the energy and potential as in corollary 7 of Ref.
\cite{aamo1}. Especially the $f_I(x,t), f_R(x,t)$ used in lemma
(\ref{lm:kers}) are also defined there.
\end{proof}

\section{Neumann Boundary Value Problem}\label{sec3}
\setcounter {equation} {0}
\setcounter {theorem} {0}
\setcounter {lemma} {0}
To stabilize the Neumann boundary value problem
\bea
\rho_t&=&a_R\rho_{xx}+b_R(x,t)\rho-a_I\iota_{xx}-b_I(x,t)\iota,\label{eq:neure}\\
\iota_t&=&a_R\iota_{xx}+b_R(x,t)\iota+a_I\rho_{xx}+b_I(x,t)\rho\label{eq:neuie}\eea
for $x\in(0,1), t\in(0,1)$, with boundary conditions
\bea \rho_x(0,t)&=&0,~\iota_x(0,t)=0\label{eq:neulbc}\\
\rho(1,t)&=&p_R(t),~\iota(1,t)=p_I(t)~,\label{eq:neurbc}\eea
the transform as (\ref{eq:rhotrs})-(\ref{eq:iotatrs}) with kernel $k(x,y,t)$ and $k_c(x,y,t)$ in the following
system
\bea
k_{xx}&=&k_{yy}+\beta(x,y,t)k+\beta_c(x,y,t)k_c+p_1k_t+p_2k_{c,t}\label{eq:neukerk}\\
k_{c,xx}&=&k_{c,yy}-\beta_c(x,y,t)k+\beta(x,y,t)k_c+q_1k_{c,t}+q_2k_t\label{eq:neukerkc}\eea
for $(x,y)\in \Omega$ with boundary conditions
\bea
k(x,x,t)&=&-\frac{1}{2}\int_0^x\beta(\gamma,\gamma,t)d\gamma\label{eq:neukintcon}\\
k_c(x,x,t)&=&\frac{1}{2}\int_0^x\beta_c(\gamma,\gamma,t)d\gamma\label{eq:neukcintcon}\\
k_x(x,0,t)&=&k_{c,x}(x,0,t)=0\\
k(0,0,t)&=&k_c(0,0,t)=0\label{eq:neugoodb}\eea where
$\beta,\beta_c,p_1,p_2,q_1,q_2,f_I$ and $f_R$ are the same as in
lemma (\ref{lm:kers}), will lead to the same system as
(\ref{eq:trrho})-(\ref{eq:bcnew}). If we can show that the above
kernels uniquely exist, all the other lemmas and theorem in section
\ref{sec2} can be applied in this case. The main result are
established as follows, in which the unique existence of the new
kernels are explicitly shown.

\begin{lemma}\label{lm:neukers}
If $a_Rb_R(x,t)+a_Ib_I(x,t)$ is analytic in $t$, the partial differential
equation system (\ref{eq:neukerk})-(\ref{eq:neugoodb}) has a unique solution satisfying
\bea
|k(x,y,t)|&\leq& Me^{M(x^2-y^2)}\label{eq:neukcon}\\
|k_c(x,y,t)|&\leq& Me^{M(x^2-y^2)},\label{eq:neukccon}\eea where $M$ is a positive
constant for $t\in (0,t_0)$ we are concerning.
\end{lemma}
\begin{proof}
The proof of this lemma is similar to that of lemma (\ref{lm:kers})
except we should find a $G(\xi, \xi, t)$ as (\ref{eq:impbc}) in
lemma (\ref{lm:kers}) since this must be used in order to obtain
$G(\xi,\eta,t)$ as (\ref{eq:geq}). Using the same substitutions as
in (\ref{eq:xietas})-(\ref{eq:deltas}), system
(\ref{eq:neukintcon})-(\ref{eq:neukcintcon}) can be transformed into
\bea G_{\xi\eta}&=&\(\delta(\xi,\eta,t)+p_1\frac{\pa}{\pa
t}\)G(\xi,\eta,t)
+\(\delta_c(\xi,\eta,t)+p_2\frac{\pa}{\pa t}\)G_c(\xi,\eta,t)\label{eq:neug}\\
G_{c,\xi\eta}&=&\(\delta(\xi,\eta,t)+q_1\frac{\pa}{\pa t}\)G_c(\xi,\eta,t)
+\(-\delta_c(\xi,\eta,t)+q_2\frac{\pa}{\pa t}\)G(\xi,\eta,t)\eea
with boundary conditions
\bea
G(\xi,0,t)&=&-\frac{1}{2}\int_0^{\xi}\delta(\tau,0,t)d\tau,\label{eq:gbc1}\\
G_c(\xi,0,t)&=&\frac{1}{2}\int_0^{\xi}\delta_c(\tau,0,t)d\tau,\\
G_{\xi}(\xi,\xi,t)&=&G_{\eta}(\xi,\xi,t),\label{eq:eqbc}\\
G_{c,\xi}(\xi,\xi,t)&=&G_{c,\eta}(\xi,\xi,t)\\
G(0,0,t)&=&G_c(0,0,t)=0.\label{eq:goodbc}\eea
Differentiating equation (\ref{eq:gbc1}) with respect to $\xi$ gives
\be G_{\xi}(\xi,0,t)=-\frac{1}{2}\delta(\xi,0,t).\label{eq:forb}\ee
Integrating equation (\ref{eq:neug}) with respect to $\eta$ from $0$ to $\xi$
and using (\ref{eq:forb}) gives
\splt{
&G_{\xi}(\xi,\xi,t)\\
=&G_{\xi}(\xi,0,t)+\int_0^{\xi}\left[\(\delta(\xi,s,t)+p_1\frac{\pa}{\pa t}\)G(\xi,s,t)
+\(\delta_c(\xi,s,t)+p_2\frac{\pa}{\pa t}\)G_c(\xi,s,t)\right]ds\\
=&-\frac{1}{2}\delta(\xi,0,t)+\int_0^{\xi}\left[\(\delta(\xi,s,t)+p_1\frac{\pa}{\pa t}\)G(\xi,s,t)
+\(\delta_c(\xi,s,t)+p_2\frac{\pa}{\pa t}\)G_c(\xi,s,t)\right]ds.}
Thus using (\ref{eq:eqbc}), we have
\splt{
&\frac{\pa G(\xi,\xi,t)}{\pa \xi}\\
=&G_{\xi}(\xi,\xi,t)+G_{\eta}(\xi,\xi,t)=2G_{\xi}(\xi,\xi,t)\\
=&-\delta(\xi,0,t)+2\int_0^{\xi}\left[\(\delta(\xi,s,t)+p_1\frac{\pa}{\pa
t}\)G(\xi,s,t) +\(\delta_c(\xi,s,t)+p_2\frac{\pa}{\pa
t}\)G_c(\xi,s,t)\right]ds.} Integrating the above equation from $0$
to $\xi$ and using (\ref{eq:goodbc}) gives the function we wanted
\splt{
&G(\xi,\xi,t)\\
=&-\int_0^{\xi}\delta(\tau,0,t)d\tau\\
&+2\int_0^{\xi}\!\!\int_0^{\tau}\left[\(\delta(\tau,s,t)+p_1\frac{\pa}{\pa
t}\)G(\tau,s,t) +\(\delta_c(\tau,s,t)+p_2\frac{\pa}{\pa
t}\)G_c(\tau,s,t)\right]dsd\tau.\label{eq:impbctwo}} Using
(\ref{eq:impbctwo}) and integrating twice (\ref{eq:neug}) first with
respect to $\eta$ from $0$ to $\eta$ and second with respect to
$\xi$ from $\eta$ to $\xi$, we can have the following integral
equation \splt{
&G(\xi,\eta,t)\\
=&-\int_0^{\eta}\delta(\tau,0,t)d\tau\\
&+2\int_0^{\eta}\!\!\int_0^{\tau}\left[\(\delta(\tau,s,t)+p_1\frac{\pa}{\pa t}\)G(\tau,s,t)
+\(\delta_c(\tau,s,t)+p_2\frac{\pa}{\pa t}\)G_c(\tau,s,t)\right]dsd\tau\\
&+\int_{\eta}^{\xi}\!\!\int_0^{\eta}\left[\(\delta(\tau,s,t)+p_1\frac{\pa
}{\pa t}\)G(\tau,s,t)+\(\delta_c(\tau,s,t)+p_2\frac{\pa
}{\pa t}\)G_c(\tau,s,t)\right]dsd\tau.}
Similarly, for $G_c(\xi,\eta,t)$ we will have
\splt{
&G_c(\xi,\eta,t)\\
=&\int_0^{\eta}\delta(\tau,0,t)d\tau\\
&+2\int_0^{\eta}\!\!\int_0^{\tau}\left[\(\delta(\tau,s,t)+q_1\frac{\pa}{\pa
t}\)G_c(\tau,s,t)
+\(-\delta_c(\tau,s,t)+q_2\frac{\pa}{\pa t}\)G(\tau,s,t)\right]dsd\tau\\
&+\int_{\eta}^{\xi}\!\!\int_0^{\eta}\left[\(\delta(\tau,s,t)+q_1\frac{\pa
}{\pa t}\)G_c(\tau,s,t)+\(-\delta_c(\tau,s,t)+q_2\frac{\pa }{\pa
t}\)G(\tau,s,t)\right]dsd\tau.} Similarly as in lemma
(\ref{lm:kers}), the above integral equations can also be shown to
have uniquely and uniformly convergent solution satisfying
(\ref{eq:neukcon}) and (\ref{eq:neukccon}).
\end{proof}

\begin{lemma}\label{lm:neuconinv}
Let $k(x,y,t)$ and $k_c(x,y,t)$ be the solution of
(\ref{eq:neukerk})-(\ref{eq:neukerkc}) and define a pair of linear bounded operator
$K$ and $K_c$ by
\bea
\trho(x,t)&=&\rho(x,t)-\int_0^x[k(x,y,t)\rho(y,t)+k_c(x,y,t)\iota(y,t)]dy\label{eq:neurhotrs}\\
\tiota(x,t)&=&\iota(x,t)-\int_0^x[-k_c(x,y,t)\rho(y,t)+k(x,y,t)\iota(y,t)]dy,\label{eq:neuiotatrs}\eea
Then,\\
1.$k$ and $k_c$ convert system (\ref{eq:neure})-(\ref{eq:neurbc}) to system
\bea
\trho_t&=&a_R\trho_{xx}+f_R(x,t)\trho-a_I\tiota_{xx}-f_I(x,t)\tiota,\label{eq:neutrrho}\\
\tiota_t&=&a_I\trho_{xx}+f_I(x,t)\trho+a_R\tiota_{xx}+f_R(x,t)\tiota,\label{eq:neutriota}\eea
for $x\in (0,1)$, with boundary conditions
\be \trho(0,t)=0,\tiota(0,t)=0,\trho(1,t)=0,\tiota(1,t)=0.\label{eq:neubcnew}\ee
2. Both $K$ and $K_c$ have linear bounded inverses.
\end{lemma}
\begin{proof}
The proof is similar to the proof of lemma (\ref{lm:coninv}).
\end{proof}

\begin{theorem}\label{thm:neu}
There exist feedback kernel $k(1,\cdot,t),~k_c(1,\cdot,t)\in
C^{2,\infty}(0,1)\times(0,t_0), t_0<1$, such that for arbitrary initial data
$\rho^0(x),~\iota^0(x)\in H^1(0,1)$, system (\ref{eq:neure})-(\ref{eq:neurbc})
with (\ref{eq:inpi})-(\ref{eq:inpr}) has a unique solution $\rho,~\iota\in
C^{2,\infty}(0,1)\times(0,t_0)$ that is exponentially stable in the
$L^2(0,1)$ and $H^1(0,1)$ norms.
\end{theorem}
\begin{proof}
The proof is similar to the proof of theorem (\ref{thm:dic}).
\end{proof}

\section{Remarks}
\setcounter {equation} {0} \setcounter {theorem} {0} \setcounter
{lemma} {0} A much more challenging but also more important task is
to stabilize the nonlinear heat equation with a
$\tilde{u}(\tilde{x},\tilde{t})$ dependent coefficient \be \frac{\pa
\tilde{u}(\tilde{x},\tilde{t})}{\pa \tilde{t}}=a_1\frac{\pa^2
\tilde{u}(\tilde{x},\tilde{t})}{\pa
\tilde{x}^2}+a_2(\tilde{x},\tilde{t},\tilde{u}(\tilde{x},\tilde{t}))\tilde{u}(\tilde{x},\tilde{t}),
\ee which models virous diffusion processes with a self-dependent
source for the diffused matter by using the boundary feedback
control method. It could be interesting to study the applicability
of the {\it method of dominant} in this case.

\section*{Acknowledgments}
The author thanks Dr. Liu Weijiu for his suggestions to consider
this problem and his valuable discussions. This work was done during
January 2005 and April 2005. The author is grateful to Dr. Vladimir
A. Miranskyy for his support during this work.

\end{document}